\documentclass[showkeys]{article}
\usepackage{amsmath,amsthm, amssymb,latexsym,amscd}
\usepackage{amsfonts}

\title{Higher order Schwarzian derivatives in interval dynamics}

\author{O.\ Kozlovski\\
\small{Mathematics Institute}\\
\small{University of Warwick}\\
\small{Coventry CV4 7AL}\\
\small{United Kingdom}
\and
D.\ Sands
\\
\small{CNRS, D\'epartement\ de Math\'ematiques}\\
\small{Universit\'e Paris-Sud}\\
\small{91405 Orsay Cedex}\\
\small{France}
}

\DeclareMathOperator{\BigO}{O}
\DeclareMathOperator{\littleo}{o}

\newtheorem{theorem}{Theorem}
\newtheorem{lemma}{Lemma}
\newtheorem{corollary}{Corollary}
\newtheorem{fact}{Fact}
\newtheorem{proposition}{Proposition}
\newtheorem{definition}{Definition}

\newcommand{\Pade}[3]{[{#1}]_{#2}^{#3}}

\newenvironment{piecewise}{\large\left\{\begin{array}{ll}}{\end{array}\right.}

\newcommand\R{{\mathbb R}}
\newcommand\C{{\mathbb C}}

\newcommand{\W}{{\cal W}}
\newcommand{\V}{{\cal V}}

\newcommand{\UHP}{{\mathbb H}}
\newcommand{\D}{{\mathbb D}}

\newcommand{\PoSch}{P}

\newcommand{\interior}{\mathrm{int}\, }

\newcommand{\degree}{\mathrm{deg}\, }
\newcommand{\dist}{\mathrm{dist}}
\newcommand{\Pock}[1]{\mathcal{P}_{#1}}
\newcommand{\Pick}[2]{\Pock{#1}({#2})}

\begin{document}

\maketitle
\abstract{%
We introduce an infinite sequence of higher order Schwarzian
derivatives closely related to the theory of monotone matrix functions.
We generalize the classical Koebe lemma to maps with positive
Schwarzian derivatives up to some order, obtaining control over
derivatives of high order.
For a large class of multimodal interval maps
we show that all inverse branches of first return maps
to sufficiently
small neighbourhoods of critical values have their
higher order Schwarzian derivatives positive up to
any given order.
}

\begin{center}
\it Dedicated to the memory of Adrien Douady
\end{center}

\section{Introduction}

Many results on interval map dynamics were first proved supposing that
the map in question has negative Schwarzian derivative.  This is a convexity
condition and as such globally constrains the possible shape of the map.
Since iterates also have negative Schwarzian derivative, estimates on the
distortion of high iterates are often greatly simplified by this assumption.
Having negative Schwarzian does not however give good control over
derivatives of high order.

The great historical failing of the negative Schwarzian theory was that many
interesting and otherwise well-behaved interval maps simply do not have negative
Schwarzian derivative everywhere.
This flaw was rectified by the discovery that a large class of interval maps are,
under some mild hypotheses, real-analytically conjugate to maps with negative
Schwarzian derivative everywhere~\cite{GS}.
An early manifestation of this was the proof~\cite{KO} that first return maps to
small neighbourhoods of critical values have negative Schwarzian derivative.

In this paper we introduce an infinite sequence of higher order Schwarzian
derivatives and prove that inverse branches of first return maps to sufficiently
small neighbourhoods of critical values have positive Schwarzian derivatives up
to some given finite order.
(For the classical Schwarzian derivative one can equally well consider branches
having negative Schwarzian derivative or inverse branches having positive
Schwarzian derivative: these are equivalent.  For higher order Schwarzian
derivatives this symmetry
disappears, and the natural concept turns out to be positivity).
We also extend the celebrated real Koebe lemma to maps with all Schwarzian
derivatives positive up to some order, obtaining control over the distortion
of derivatives of high order.

Let $d$ be a positive integer and $f$ a map which is $2d+1$ times differentiable at $x$.
Let $\mathcal R$ be the rational map of degree at most $d$
that coincides with $f$ to order $2d$ at $x$, i.e.\ for which
${\mathcal R}(x) = f(x)$, $D{\mathcal R}(x) = Df(x)$, \ldots,
$D^{2d}{\mathcal R}(x) = D^{2d}f(x)$.
Such a rational map may not exist, but if it does exist then it is unique.
It is called the $d$'th (diagonal) Pad\'e approximant~\cite{PA} to $f$ at $x$.
The Schwarzian derivative of $f$ at $x$ of order $d$ is then defined to be
$S_d(f)(x) \equiv D^{2d+1}({\mathcal R}^{-1} \circ f)(x)$.
Only a local inverse of $\mathcal R$ being needed,
this makes sense as long as $Df(x) \neq 0$.
It is well-known, and easily checked, that this definition gives the
classical Schwarzian derivative
$S_1(f)(x) = \frac{D^3 f(x)}{Df(x)} -
\frac{3}{2}\left(\frac{D^2 f(x)}{Df(x)}\right)^2$
when $d=1$.
Define $S_0(f) \equiv 1$ for convenience.

In order to simplify the exposition, existence of higher order
Schwarzian derivatives will be indicated implicitly with the
convention that an expression such as $S_d(f)(x) < 0$ is short-hand for:
$S_d(f)(x)$ exists and $S_d(f)(x) < 0$.

It is essential to the theory of the classical Schwarzian derivative
that iterates of maps with negative Schwarzian derivative also have
negative Schwarzian derivative, as follows from the composition formula
$S_1(f \circ g) = S_1(f) \circ g \, (Dg)^2 + S_1(g)$.  There is also
a composition formula for the higher order Schwarzian derivative $S_d$
(lemma~\ref{lcomp}), but it contains an extra term coming from the fact
that the set of rational maps of degree $d > 1$ is not closed under
composition.
This term disappears when composing with a M\"obius transformation
$M$ (i.e.\ a rational map of degree $1$): post-composition has no effect
while pre-composition yields
$S_d(f \circ M) = S_d(f) \circ M \, (DM)^{2d}$.
In the special case of maps with all Schwarzian derivatives of
order less than $d$ non-negative, the extra term is non-negative
(proposition~\ref{pcomp}).
This makes estimating the $d$'th order Schwarzian derivative of a
long composition of such functions feasible.
It also shows that this class of maps is closed under composition.
The ultimate origin of this is the fact
(corollary~\ref{cschpick}) that
$S_1(f)(x) > 0, \ldots, S_{d-1}(f)(x) > 0$ if and only if
the $d$'th Pad\'e approximant to $f$
at $x$ exists, has degree $d$,
and maps the complex upper half-plane into itself
(if $Df(x) > 0$) or into the complex lower half-plane (if
$Df(x) < 0$).

The following result explains why inverse branches of many well-known
one-dimensional maps (logistic maps for example) have Schwarzian
derivatives of all orders positive:
they map the complex upper half-plane into itself.  This is
almost but not quite the same as being in the Epstein class.
Given an open interval $U$, let $\PoSch_d(U)$ consist of those
functions $f : U \to \R$ with $2d+1$ derivatives,
$Df \neq 0$ and $S_1(f) \geq 0$, $\ldots$, $S_d(f) \geq 0$
everywhere;
set $\PoSch_\infty(U) = \cap_{d=1}^\infty \PoSch_d(U)$.

\begin{proposition} \label{puni}
Let $\phi : U \to \R$ be an increasing $C^\infty$ diffeomorphism
onto its image, where $U$ is an open real interval.
Then $\phi \in \PoSch_\infty(U)$ if and only if $\phi$ extends to a
holomorphic map $\phi : \C \setminus (\R \setminus U) \to \C$ which
maps the complex upper half-plane into itself (Pick class).
\end{proposition}

It has often been observed that smooth maps become ``increasingly
holomorphic'' when iterated.  In the light of the preceding
proposition, one way of formalising this observation is to say
that Schwarzian derivatives of ever increasing order
become positive for inverse branches under iteration.
Our main result shows that this is indeed the case near critical values:

\begin{theorem} \label{tmain} Let $f : I \to I$ be a $C^{2d+1}$ map of a
  non-trivial compact interval, and let all critical points of $f$ be non flat.
  Then for any critical point $c$ of $f$ which is not in the basin of a
  periodic attractor there exists a neighbourhood $X$ of $c$ such that if
  $f^s(x) \in X$, for some $x \in I$ and $s\ge 0$ with
  $Df^{s+1}(x)\neq 0$, then $S_k(f^{-(s+1)})(f^{s+1}(x)) > 0$ for all
  $k=1,\ldots,d$ (local inverse near $x$).
\end{theorem}

Recall that a critical point $c$ is said to be \emph{non flat}
if $f$ can be decomposed near $c$ as $f=\psi\circ P \circ \phi$
where $\phi$ (resp. $\psi$) is a 
$C^{2d+1}$ diffeomorphism from a neighbourhood of $c$ (resp. $0$)
onto a neighbourhood of $0$ (resp. $f(c)$) and $|P(x)|=|x|^\alpha$
for some $\alpha > 1$ and all small $x$.

We have so far been unable to prove, when $d > 1$, the global result
corresponding to theorem~\ref{tmain}, namely
that if all periodic points are hyperbolic repelling then $f$ can
be real-analytically conjugated to a map with all
Schwarzian derivatives of order $d$ or less positive
everywhere for all inverse branches
(the conjugacy would depend on $d$).

The celebrated Koebe lemma for univalent maps $f : U \to V$,
where $U$ and $V$ are simply connected domains in $\C$,
states that $f$ has bounded distortion on any simply connected domain
$A$ which is compactly contained in $U$:
if $x \in A$ and $y \in A$ then $|Df(x)|/|Df(y)|$ is bounded by a
constant depending only on the modulus of $U\setminus A$.
Higher derivatives are also controlled:
after appropriately normalising the domain $U$,
the ratio $|D^kf(x)|/|Df(x)|$ is
again bounded by a constant depending only on $k$ and the modulus of
$U\setminus A$.

A real counterpart to the complex Koebe lemma exists but lacks the
control of higher order derivatives.
Here we show that the real Koebe lemma can be naturally generalised,
and control of higher derivatives achieved, in a similar way to the complex
Koebe lemma:

\begin{theorem} \label{koebe}
  Let $d$ be a positive integer, $U$ an open interval and
  $m$ and $n$ integers such that $n$ is odd and $1 \leq n \leq m \leq 2d$.
  If $f \in \PoSch_d(U)$ then
  \begin{equation} \label{kbound}
  |D^mf(x)| \leq \frac {n!}{m!}
  \,\dist(x, \partial U)^{-(m-n)} \, |D^nf(x)|
  \end{equation}
  for every $x$ in $U$,
  where $\dist(x, \partial U)$ is the distance from $x$ to the boundary of $U$.
  The constant in (\ref{kbound}) is exact and is achieved on a M\"{o}bius
  transformation.
\end{theorem}

Note that the inequality in the theorem has the right scaling
properties with respect to affine coordinate changes.

Many additional useful facts about maps in $\PoSch_d$ can be
deduced from the theory of monotone matrix functions~\cite{DO}.
This is because maps in $\PoSch_d$
with positive derivative turn out (lemma~\ref{lmm}) to be
exactly the monotone matrix functions of order $d+1$.
Consider for example theorem VII.V from~\cite{DO}.  Dividing the
rows and columns of the Pick matrix appropriately turns the
matrix elements into cross-ratios;
cross-ratios being unaffected by sign changes, the resulting
matrix is positive for any function in $\PoSch_d$, not just
those that are increasing.  The result is a generalisation of
the well-known cross-ratio contraction property of maps
with positive Schwarzian derivative:

\begin{proposition}
Let $d$ be a positive integer, $U$ an open interval and take
distinct points $\lambda_1, \ldots, \lambda_{d+1}$ in $U$.
If $f \in \PoSch_d(U)$ then all eigenvalues of the matrix
\begin{equation}
\left[
\sqrt{
\frac{\left(\frac{f(\lambda_i) - f(\lambda_j)}{\lambda_i - \lambda_j}\right)^2}
{Df(\lambda_i) Df(\lambda_j)}
}
\right]
\end{equation}
(diagonal elements are equal to $1$)
are non-negative: the matrix is positive.
\end{proposition}

Some other properties likely to be useful for dynamics can be
found in~theorem VII.II, chapter XIV and theorem II.I of~\cite{DO}.
\\
\\
\emph{Organisation.}
Section~\ref{ratapprox} sketches the essentials of the
theory of Pad\'e approximation as used in this paper and
describes some elementary properties of higher order
Schwarzian derivatives.
Section~\ref{spick} introduces the Pick algorithm,
a degree reduction technique useful for proving
results by induction on the order of the Schwarzian
derivative.
This technique is used in section~\ref{sratpick} to
characterise the real rational maps that preserve the
complex upper half plane in terms of their higher order
Schwarzian derivatives at a point.
At this point the theory is sufficiently developed that
obtaining an effective composition formula is trivial ---
this is done in section~\ref{sec:composition-formula}.
Another application is the proof in section~\ref{smonmat}
that the increasing functions in $\PoSch_d$ are exactly
the monotone matrix functions of order $d+1$.
Proposition~\ref{puni} is thus a restatement of Loewner's
theorem~\cite{DO}.
Section~\ref{skoebe} uses the theory of monotone matrix functions
to deduce the generalised Koebe lemma from the integral representation
for Pick functions.
Finally, in section~\ref{smain}, we prove the main theorem
using the a priori bounds of~\cite{SV} and the
results from~\cite{ST} on Epstein class approximation.

\section{Rational Approximation} \label{ratapprox}

This section contains a quick introduction to the classical
theory of rational approximation.  Most of the results are
reformulations of well-known properties.

A \emph{rational map} is a fraction ${\mathcal R} = \frac{p}{q}$
where $p$ and $q$ are polynomials and $q$ is not identically zero.
We consider two rational maps ${\mathcal R}_1 = \frac{p_1}{q_1}$
and ${\mathcal R}_2 = \frac{p_2}{q_2}$ to be equal if the polynomials
$p_1 q_2$ and $p_2 q_1$ are equal.  This means that common polynomial
factors in the numerator and the denominator can be cancelled
without changing the rational map.  When viewing $\mathcal R$ as a
function we will suppose that $p$ and $q$ are relatively prime, i.e.\
that $z \mapsto p(z)/q(z)$ has no removable singularities on the
Riemann sphere.  The \emph{degree} of $\mathcal R$, denoted
$\degree \mathcal R$, is the maximum of the degrees of $p$ and $q$
when $p$ and $q$ are relatively prime.

Two functions $f$ and $g$ are said to \emph{coincide to order $N$} at some
point $x$ if $f(x) = g(x)$, $Df(x) = Dg(x)$, $\ldots$, $D^Nf(x) = D^Ng(x)$.

\begin{lemma}[Uniqueness] \label{lwd1}
Let ${\mathcal R}_1$ (resp.\ ${\mathcal R}_2$) be a rational map of
degree at most $d_1$ (resp.\ $d_2$).
If  ${\mathcal R}_1$ and  ${\mathcal R}_2$ are finite at some point $x$,
and coincide to order $d_1+d_2$ there, then ${\mathcal R}_1 = {\mathcal R}_2$.
\end{lemma}
\begin{proof}
See \cite[Theorem 1.1]{BA}.
Without loss of generality $x = 0$.
Let ${\mathcal R}_i = p_i/q_i$ where $p_i$, $q_i$ are polynomials
of degree at most $d_i$, $i=1,2$.
The hypothesis that ${\mathcal R}_1$ and ${\mathcal R}_2$ coincide
to order $d_1+d_2$ is equivalent to
$p_1(z)/q_1(z) - p_2(z)/q_2(z) = \BigO(z^{d_1+d_2+1})$.
Multiplying through by the denominators gives
$p_1(z)q_2(z)-p_2(z)q_1(z) = \BigO(z^{d_1+d_2+1})$.
The left-hand side is a polynomial of degree at most $d_1+d_2$
so must in fact be identically zero: $p_1 q_2 = p_2 q_1$.
\end{proof}

Let $d$ be a non-negative integer and $f$ a map which is $2d$ times
differentiable at $x$.  Recall that the $d$'th (diagonal) Pad\'e
approximant to $f$ at $x$, denoted $\Pade{f}{x}{d}$, is the rational
map of degree at most $d$ that coincides with $f$ to order $2d$ at $x$,
if such a rational map exists.  By lemma~\ref{lwd1} there is at most
one such a rational map.

The classical sufficient condition for the existence of the $d$'th
Pad\'e approximant is the non-vanishing of a certain Hankel determinant.
The empty determinant is considered equal to $1$.

\begin{lemma}[Existence] \label{lwd3}
Let $d$ be a non-negative integer and $x$ some point.
Given numbers $f_0$, $f_1$, $\ldots$, $f_{2d}$ there exists a
rational map ${\mathcal R}$ of degree exactly $d$ such that
${\mathcal R}(x) = f_0$,
$D{\mathcal R}(x) = f_1$, $\ldots$,
$D^{2d}{\mathcal R}(x) = f_{2d}$ if and only if
\begin{equation} \label{enormal}
\det\left[ \begin{array}{cccc}
\frac{f_1}{1!} & \frac{f_2}{2!}         & \cdots & \frac{f_d}{d!}           \\
\frac{f_2}{2!} & \frac{f_3}{3!}         & \cdots & \frac{f_{d+1}}{(d+1)!}   \\
\vdots         &                        &        & \vdots                   \\
\frac{f_d}{d!} & \frac{f_{d+1}}{(d+1)!} & \cdots & \frac{f_{2d-1}}{(2d-1)!}
\end{array} \right]
\neq 0.
\end{equation}
\end{lemma}
\begin{proof}
See \cite[Theorem 2.3, points (1) and (5)]{BA}
or \cite[Proposition 3.2]{NS}.
Without loss of generality $x = 0$.
Write $F_k = f_k/k!$ and let
$f(z) = F_0 + F_1 z + \cdots F_{2d} z^{2d} + \BigO(z^{2d+1})$.
The existence of an appropriate $\mathcal R$, but of degree at
most (rather than exactly) $d$,
is equivalent to the existence of a pair $(p,q)$ of polynomials
such that
\begin{eqnarray}
\mbox{degree } p \leq d &,& \mbox{degree } q \leq d \nonumber \\
q(0) & = & 1 \label{ematch} \\
F(z)q(z) - p(z) & = & \BigO(z^{2d+1}). \nonumber
\end{eqnarray}
Writing $q(z) = 1 + Q_1 z + \cdots Q_d z^d$,
$p(z) = P_0 + P_1 z + \cdots P_d z^d$ and expressing
equation~\ref{ematch} in terms of coefficients yields
\begin{eqnarray}
0 = F_l + \sum_{i=1}^d Q_i F_{l-i}, & l = d+1, \ldots, 2d, \label{A.1} \\
P_l = F_l + \sum_{i=1}^l Q_i F_{l-i}, & l = 0, \ldots, d. \label{A.2}
\end{eqnarray}
The key is equation~\ref{A.1}, since if it can be solved for
$Q_1, \ldots, Q_d$ then the values of $P_0, \ldots, P_d$
are determined by (\ref{A.2}).
Thus equation~\ref{ematch} has a unique solution if and only if
the linear part of (\ref{A.1}) is invertible.
But this is the case if and only if it has non-zero determinant,
which is exactly~(\ref{enormal}).

The preliminaries now being in place, let us prove the result.

Suppose (\ref{enormal}) holds; then
(\ref{ematch}) has a unique solution $(p,q)$.
Uniqueness implies that at least one of $p$, $q$ has degree $d$,
since if not the pair $(p(z)(1+z), q(z)(1+z))$ would be another
solution to (\ref{ematch}).  Uniqueness also implies that
$p$ and $q$ are relatively prime since otherwise another solution
could be obtained by cancelling common factors
(this reasoning uses $q(0) \neq 0$).

Conversely, suppose $\mathcal R$ exists with degree exactly $d$.
Then there exists a solution $(p,q)$ of
(\ref{ematch}) where $p$ and $q$ are relatively prime polynomials,
and one of $p,q$ has degree $d$.  Suppose there exists another
solution $(\tilde{p}, \tilde{q})$ of (\ref{ematch}).  Since
$\tilde{p} q = p \tilde{q}$ and $p$ and $q$ are relatively
prime, there must exist some non-zero polynomial $r$ such
that $\tilde{p} = r p$, $\tilde{q} = r q$.  Now
$\tilde{p}$ and $\tilde{q}$ have degree at most $d$
while one of $p, q$ has degree $d$, so $r$ must have
degree zero, i.e.\ be a constant.
From $q(0)=1=\tilde{q}(0)$ it follows that the constant $r$
is in fact $1$, i.e.\ $\tilde{p}=p$ and $\tilde{q}=q$.
In other words (\ref{ematch}) has a unique solution,
which means that (\ref{enormal}) holds.
\end{proof}

Write
\begin{equation}
M_d(x,f) \equiv
\left[ \begin{array}{cccc}
\frac{Df(x)}{1!}   & \cdots & \frac{D^df(x)}{d!} \\
\vdots             &        & \vdots             \\
\frac{D^df(x)}{d!} & \cdots & \frac{D^{2d-1}f(x)}{(2d-1)!}
\end{array} \right].
\end{equation}

\begin{corollary} \label{dzero}
Let $d$ be a non-negative integer, $U$ an open interval and
$f : U \to \R$ a function with $2d+1$ derivatives.
Then $\det M_{d+1}(x,f)=0$ for all $x \in U$ if and only if
$f$ is a rational map of degree at most $d$.
\end{corollary}
\begin{proof}
Suppose $f$ is a rational map of degree at most $d$.
If $\det M_{d+1}(x,f) \neq 0$ then by lemma~\ref{lwd3} there
is a rational map $\mathcal R$ of degree $d+1$ which
coincides with $f$ to order $2d+2$ at $x$.
Then lemma~\ref{lwd1} shows that $f = {\mathcal R}$,
which is impossible because their degrees differ.
Thus $\det M_{d+1}(x,f) =0$ for all $x \in U$.
For the converse, the difficulty is that the differential
equation to be solved is singular at points where $\det M_d(x,f)=0$.

\emph{Claim}.  Suppose $\det M_d(x,f)\neq 0$, $\det M_{d+1}(x,f) = 0$
for all $x$ in some open interval $V \subseteq U$.
Then $f$ is a rational map of degree $d$ on $V$.

\emph{Proof of claim}.
Observe that $\det M_{d+1}(x,f) = 0$ is an ordinary differential
equation for $f$ with highest term $D^{2d+1}(f)(x)$.  The
coefficient of this term is $\det M_d(x,f)$ which is
non-zero by hypothesis, so the differential equation is
non-singular: there is local existence and uniqueness.
Let $\mathcal R$ be the Pad\'e approximant to $f$ of order
$d$ at some point $p$ of $V$
($\mathcal R$ exists and has degree $d$ by lemma~\ref{lwd3}).
Note that $\det M_{d+1}(x,{\mathcal R}) =0$ for all $x \in V$,
i.e. $\mathcal R$ is a solution of the differential equation
(this was proved above).
By definition $\mathcal R$ coincides with $f$ to order $2d$ at
$p$, so $f$ and $\mathcal R$ have the same initial values at $p$
and thus coincide throughout $V$.  This proves the claim.

So suppose $\det M_{d+1}(x,f)=0$ for all $x \in U$.
First consider the case when $M_d(x,f) \neq 0$ for some point
$x \in U$ (this is always the case if $d=0$).
Then $M_d(\cdot,f) \neq 0$ at every point of $U$.  Indeed,
let $V$ be the maximal open interval around $x$ on which
$M_d(\cdot,f) \neq 0$.
By the claim, $f$ coincides with a rational map $\mathcal R$ of degree
$d$ on $V$.  If $V \neq U$ then there is some boundary
point $p$ of $V$ in $U$.  By continuity,
$D^k({\mathcal R})(p) = D^k(f)(p)$ for $0 \leq k \leq 2d$,
so $\det M_d(p,f)\neq 0$ by lemma~\ref{lwd3}.
This implies that $p \in V$, a contradiction since $V$ is open.
Thus $V=U$ and $f$ is a rational map of degree $d$.

Now suppose that $\det M_{d}(x,f) = 0$ for all $x \in U$.  This
reduces $d$ by $1$ in the hypotheses.
By induction $f$ is a rational map of degree at most $d-1$.
\end{proof}

\begin{definition}
We say that $f$ is \emph{normal} of order $d$ at $x$ if the $d$'th
Pad\'e approximant to $f$ at $x$ exists and has degree exactly $d$.
\end{definition}

\begin{corollary} \label{cnorm}
Let $f$ have $2d$ derivatives at $x$.  Then $f$ is normal of
order $d$ at $x$ if and only if $\det M_d(x,f) \neq 0$.
\end{corollary}

If $f$ is not normal of order $d$ then the $d$'th Pad\'e
approximant may nonetheless exist, but if so it is simply
equal to the $d-1$'st Pad\'e approximation:

\begin{corollary} \label{cnormal}
Let $f$ have $2d+2$ derivatives at $x$.  Suppose the
$d+1$'st Pad\'e approximant $\Pade{f}{x}{d+1}$ exists.
If $f$ is not normal of order $d+1$ at $x$ then the
$d$'th Pad\'e approximant $\Pade{f}{x}{d}$ exists
and $\Pade{f}{x}{d+1} = \Pade{f}{x}{d}$.
If $f$ is normal of order $d+1$ at $x$ then either $f$ is
normal of order $d$ at $x$ or the $d$'th Pad\'e approximant
to $f$ at $x$ does not exist.
\end{corollary}
\begin{proof}
See also~\cite[Theorem 2.3]{BA}.
If $f$ is not normal of order $d+1$ then $\Pade{f}{x}{d+1}$ has degree
at most $d$.  Since it satisfies the conditions to be the
$d$'th Pad\'e approximant, by uniqueness it is $\Pade{f}{x}{d}$.
Now suppose $f$ is normal of order $d+1$ and the $d$'th Pad\'e approximant
$\Pade{f}{x}{d}$ exists.  If $f$ is not normal of order $d$ then
$\Pade{f}{x}{d}$ has degree at most $d-1$.  Since
$\Pade{f}{x}{d}$ and $\Pade{f}{x}{d+1}$ coincide to order
$2d$ at $x$, as both coincide with $f$ to at least
that order, they are equal by lemma~\ref{lwd1}.
This contradicts $\Pade{f}{x}{d+1}$ having degree $d+1$.
\end{proof}

Another way of viewing this result is as follows:
suppose $f$ is normal of order $d$ but not of order $d+1$.
Let $N \in \{d+1, d+2, \ldots\}$ be minimal such that
the $N$'th Pad\'e approximant does not exist, or set
$N=\infty$ if Pad\'e approximants exist of all orders.
Then the %
approximants of orders $d < k < N$ are
all equal to $\Pade{f}{x}{d}$.

Recall that the Schwarzian derivative of $f$
at $x$ of order $d$ is defined to be $S_d(f)(x) =
D^{2d+1}({(\Pade{f}{x}{d})}^{-1} \circ f)(x)$.
An alternative definition is
$S_d(f)(x) = D^{2d+1}(f - \Pade{f}{x}{d})(x) / Df(x)$.
The equivalence of the two expressions is readily derived
by induction from the fact that $f$ and $\Pade{f}{x}{d}$
coincide to order $2d$ at $x$.

\begin{lemma}[Schwarzian formula] \label{lform}
Let $d$ be a non-negative integer and $x$ some point.
If $f$ has $2d+1$ derivatives at $x$, $Df(x) \neq 0$,
and $f$ is normal of order $d$ at $x$, then
\begin{equation} \label{esderiv}
S_d(f)(x) = (2d+1)! \frac{\det M_{d+1}(x,f)}{Df(x) \det M_d(x,f)}.
\end{equation}
\end{lemma}
\begin{proof}
Let ${\mathcal R} = \Pade{f}{x}{d}$ be the $d$'th Pad\'e approximant
to $f$ at $x$.  Recall the alternative definition
$S_d(f)(x) = (D^{2d+1}f(x) - D^{2d+1}{\mathcal R}(x))/Df(x)$
of the Schwarzian derivative of order $d$.
Write
\begin{equation}
A = \left[ \begin{array}{ccc}
F_1    & \cdots & F_d    \\
\vdots &        & \vdots \\
F_d    & \cdots & F_{2d-1}
\end{array} \right],
\end{equation}
\begin{equation}
C = \left[ \begin{array}{ccc}
F_{d+1} & \cdots & F_{2d}
\end{array} \right]
\end{equation}
and $D = F_{2d+1}$, where $F_k = D^k f(x)/k!$.
Let ${\mathcal R} = \frac{p}{q}$ where
$p(z) = P_0 + P_1 z + \cdots P_d z^d$ and
$q(z) = 1 + Q_1 z + \cdots Q_d z^d$
are polynomials of degree at most $d$.
Equation~\ref{A.1} can be rewritten as
$[Q_d, \cdots, Q_1]^T = - A^{-1} C^T$.
Since $D^{2d+1}{\mathcal R}(x)/(2d+1)! =
- C [Q_d, \cdots, Q_1]^T$,
it follows that
$(D^{2d+1}f(x) - D^{2d+1}{\mathcal R}(x))/(2d+1)! = D - C A^{-1} C^T$.
Applying the well-known formula
\begin{equation}
\det \left[
\begin{array}{cc}
A & B \\
C & D
\end{array}
\right] = \det A \det (D - C A^{-1} B)
\end{equation}
for the determinant of a conformally partitioned block matrix
with $B = C^T$ immediately yields
\begin{equation}
D - C A^{-1} C^T = \frac{\det M_{d+1}(x,f)}{\det M_d(x,f)}.
\end{equation}
Thus
\begin{equation}
S_d(f)(x) = \frac{D^{2d+1}f(x) - D^{2d+1}{\mathcal R}(x)}{Df(x)}
= (2d+1)! \frac{\det M_{d+1}(x,f)}{Df(x) \det M_d(x,f)}.
\end{equation}
\end{proof}

\begin{lemma} \label{lzero}
Let $d$ be a non-negative integer, $U$ an open interval and
$f : U \to \R$ a function with $2d+1$ derivatives.  Suppose
$Df(x) \neq 0$ for all $x\in U$.
Then $S_d(f)$ is identically zero on $U$ if and only if $f$ is
a rational map of degree at most $d$.
\end{lemma}
\begin{proof}
If $f$ is a rational map of degree at most $d$, then $f$ is
the $d$'th Pad\'e approximant to itself at any point,
so $S_d(f)=0$ by definition.
So suppose $S_d(f)$ is identically zero on $U$.
If $f$ is not normal of order $d$ at any point of $U$ (this does not
happen if $d=0$) then $\det M_d(x,f) = 0$ for all $x \in U$, so
$f$ is a rational map of degree at most $d-1$ by lemma~\ref{dzero}.
So suppose $f$ is normal of order $d$ at some point $x$ of $U$.  Then,
as in the proof of lemma~\ref{dzero}, $f$ is normal of order $d$ at
every point of $U$.  Indeed, let $V$ be the maximal open interval
around $x$ on which $f$ is normal of order $d$.  Then formula~\ref{esderiv}
is valid on $V$, so from $S_d(f)\equiv 0$ it follows that
$\det M_{d+1}(\cdot,f) = 0$ on $V$.
It was shown in the proof of lemma~\ref{dzero} that
$f$ then coincides on $V$ with a rational map $\mathcal R$
of degree exactly $d$.
If $V \neq U$ then there is some boundary point $p$ of $V$ in $U$.
By continuity, $f$ and $\mathcal R$ coincide to order $2d$ at $p$,
so $\mathcal R$ is the $d$'th Pad\'e approximant to $f$ at $p$.
Thus $f$ is normal of order $d$ at $p$ by definition, i.e.\ $p \in V$.
This contradicts $V$ being open.
Hence $V=U$ and $f$ is a rational map of degree $d$.
\end{proof}

\section{The Pick Algorithm} \label{spick}

Let $\D = \{z \in \C : \left|z\right| < 1\}$ be the open unit disk
in the complex plane.
The \emph{Schur class} consists of all holomorphic
functions $\D \to \overline{\D}$.
Given a holomorphic map $F : \D \to \D$, the Schur
algorithm generates a new holomorphic map
$\tilde{F} : \D \to \overline{\D}$ as follows.
Let $M_F$ be the M\"{o}bius transformation
$M_F : z \mapsto (z-F(0))/(1-\overline{F(0)}z)$.
This preserves $\D$ and maps $F(0)$ to $0$.
The function $z \mapsto M_F(F(z))/z$
has a removable singularity at $z=0$,
so extends to a holomorphic function
$\tilde{F} : \D \to \C$.
The Schwartz lemma shows that in fact
$\tilde{F} : \D \to \overline{\D}$.
If $\tilde{F}$ is not constant then
$\tilde{F} : \D \to \D$.
Applying the algorithm iteratively results
in a finite or infinite sequence of Schur maps
that terminates with a constant function if finite.

The Schur algorithm uses $0$ as a distinguished point of $\D$ and
a particular choice of M\"{o}bius transformation taking $F(0)$ to $0$.
Other choices produce different sequences of maps.
Moving the distinguished point towards the boundary of $\D$,
normalising with M\"{o}bius transformations,
and passing to the limit results in a version of the Schur
algorithm for which the distinguished point lies on the unit circle.
This is the Pick algorithm studied in this section.

Rather than work with $\D$, $\overline{\D}$ and the unit circle,
it is more convenient to use the conformally equivalent
complex upper half-plane $\UHP = \{z \in \C : \Im(z) > 0\}$,
its closure in the Riemann sphere
$\overline{\UHP} = \{z \in \C : \Im(z) \geq 0\} \cup \{\infty\}$,
and the extended real-line $\R \cup \{\infty\}$.
Section~\ref{sratpick} makes use of the Pick algorithm in the complex
plane.
Here we consider real-valued functions defined on a real neighbourhood
of a point $x \in \R$ since this suffices for our applications.

\begin{definition}
Let $f$ be twice differentiable at $x$ with $Df(x) \neq 0$.
The \emph{Pick algorithm} based at $x$ transforms $f$ into
\begin{equation}
\Pick{x}{f} : z \mapsto
\begin{piecewise}
\frac{1 - Df(x) \frac{z-x}{f(z)-f(x)}}{z-x} & z \neq x \\
\frac{D^2f(x)}{2\,Df(x)} & z = x.
\end{piecewise}
\end{equation}
\end{definition}

Note that $\Pick{x}{f}$ is continuous at $x$ (in general two
derivatives are lost at $x$).
If $f(x)$ and $Df(x)$ are known, then $f$ can be
recovered from $\tilde{f} = \Pick{x}{f}$:

\begin{definition}
Let $\tilde{f}$ be continuous at $x$,
and take some $A \in \R$ and $\mu \in \R \setminus \{0\}$.
The \emph{inverse Pick algorithm} based at $x$ transforms
$\tilde{f}$ into
\begin{equation}
f : z \mapsto A + \frac{\mu(z-x)}{1-(z-x)\tilde{f}(z)}.
\end{equation}
\end{definition}
This is indeed an inverse: $\Pick{x}{f} = \tilde{f}$.
Note that $f$ is twice differentiable at $x$,
$f(x) = A$ and $Df(x) = \mu \neq 0$ (in general two derivatives
are gained at $x$).

Clearly $f$ is a rational map if and only if $\Pick{x}{f}$ is.
If they are rational, it is straightforward to show that
$\degree f = 1 + \degree \Pick{x}{f}$.  This uses the standing
assumptions that $f$ is real (which implies that $f(x)$ is finite) and $Df(x) \neq 0$.

It can be helpful to think of the Pick algorithm in terms of
continued fractions.
Applying the inverse Pick algorithm $d$ times to $\tilde{f}$
results in the Jacobi-type continued fraction
\begin{equation*}
f : z \mapsto
A_0 + \cfrac{\mu_0(z-x)}{1 - (z-x) A_1 - \cfrac{\mu_1(z-x)^2}{1 - (z-x) A_2 - \cfrac{\mu_2(z-x)^2}{\,\ddots\,-\cfrac[l]{\mu_{d-1}(z-x)^2}{1 - (z-x)\tilde{f}(z)}}}}
\end{equation*}
\begin{equation} \label{confrac}
\quad\quad A_0, \ldots, A_{d-1} \in R, \quad \mu_0, \ldots, \mu_{d-1} \in \R \setminus \{0\}.
\end{equation}
Since we only use continued fractions to illustrate results
rather than prove them, we have felt free to state their
properties without justification.
Observe that $f$ is $2d$ times differentiable at $x$ and
$\Pock{x}^d(f) = \tilde{f}$.
More: $f$ is normal of orders $1, \ldots, d$ at $x$.
Conversely, if $f$ is $2d$ times differentiable at $x$
and $f$ is normal of orders $1, \ldots, d$ at $x$ then
$f$ can be written in the form~(\ref{confrac}) with
$\tilde{f}$ continuous at $x$.
The $k$'th convergent of equation~\ref{confrac}
(obtained by setting $\mu_k = 0$) is exactly
the $k$'th Pad\'e approximant to $f$ at $x$
($0 \leq k < d$).  These properties are the essence of:

\begin{lemma}
Let $d$ be a positive integer and $f$ a map with is
$2d$ times differentiable at $x$ with $Df(x) \neq 0$.
The $d$'th Pad\'e approximant $\Pade{f}{x}{d}$ to $f$
at $x$ exists if and only if the $d-1$'st Pad\'e
approximant $\Pade{\Pick{x}{f}}{x}{d-1}$ to $\Pick{x}{f}$
exists, and then
$\Pick{x}{\Pade{f}{x}{d}} = \Pade{\Pick{x}{f}}{x}{d-1}$.
\end{lemma}
\begin{proof}
Write $\tilde{f}$ for $\Pick{x}{f}$.
The case $d=1$ is immediate: the $0$'th Pad\'e approximant
always exists and the $1$'st Pad\'e approximant to $f$ exists
because $Df(x) \neq 0$.  The formula connecting the two
is trivial.
So suppose $d > 1$ and that $\Pade{\tilde{f}}{x}{d-1}$ exists ---
denote it by $\mathcal T$.  Define $\mathcal R$ via the inverse
Pick algorithm:
$\mathcal R(z) = f(x) + Df(x)(z-x)/(1 - (z-x) {\mathcal T}(z))$.
Note that $\degree {\mathcal R} = 1 + \degree {\mathcal T}$.
Then
\begin{equation} \label{upeq}
f(z) - {\mathcal R}(z) =
\frac{1}{Df(x)}\frac{f(z) - f(x)}{z-x}
\frac{{\mathcal R}(z) - {\mathcal R}(x)}{z-x}
(z-x)^2 (\tilde{f}(z) - {\mathcal T}(z)).
\end{equation}
By definition $\tilde{f}$ and $\mathcal T$ coincide to order $2(d-1)$,
i.e.\ $\tilde{f}(z) - {\mathcal T}(z) = \littleo((z-x)^{2(d-1)})$.
Then $f(z) - {\mathcal R}(z) = \littleo((z-x)^{2d})$
by equation~\ref{upeq}, which means that $f$ and $\mathcal R$ coincide
to order $2d$.
Thus $\mathcal R$ is the $d$'th Pad\'e approximant to $f$ at $x$,
as desired.
The case when it is the $d$'th Pad\'e approximant
to $f$ at $x$ that is initially known to exist
is left to the interested reader.
\end{proof}

\begin{corollary} \label{schsch}
Let $d$ be a positive integer and $f$ a map which is $2d+1$ times
differentiable at $x$ with $Df(x) \neq 0$.
Then
\begin{equation}
S_d(f)(x) = 2d(2d+1) D\Pick{x}{f}(x) S_{d-1}(\Pick{x}{f})(x). \label{equpsch}
\end{equation}
Included in this is that
$S_d(f)$ exists if and only if $S_{d-1}(\Pick{x}{f})$ exists.
\end{corollary}
\begin{proof}
The statement about existence is immediate from the previous
lemma.
For the formula, revisit the proof of the previous lemma.
Writing
$\tilde{f}(z) - {\mathcal T}(z) = \alpha (z-x)^{2d-1} + \littleo((z-x)^{2d-1})$,
observe that
$\alpha = D\tilde{f}(x) S_{d-1}(\tilde{f})(x) / (2d-1)!$.
This observation is precisely the alternative definition of the
higher Schwarzian derivative from section~\ref{ratapprox}.
Likewise, $f(z) - {\mathcal R}(z) =
Df(x) \alpha (z-x)^{2d+1} + \littleo((z-x)^{2d+1})$
--- which follows from equation~\ref{upeq} ---
means $S_d(f)(x) = \alpha (2d+1)!$.
This is the same as~(\ref{equpsch}).
\end{proof}

In terms of the continued fraction representation~(\ref{confrac}),
this says that
$\mu_k = \frac{S_k(f)(x)}{2k(2k+1)S_{k-1}(f)(x)}$ for $1 \leq k < d$.

\section{Rational Pick maps} \label{sratpick}

In this section we use the Pick algorithm to characterise the real
rational maps in the Pick class as those with their Schwarzian
derivatives of all orders non-negative.  Remarkably, if they are
non-negative at a single point then they are non-negative everywhere.

A rational map $\mathcal R$ is \emph{real}
if it can be written as a ratio of polynomials with only real
coefficients.  This is equivalent to $\mathcal R$ being real
or infinite valued everywhere on the real line.
By uniqueness, Pad\'e approximants to real-valued maps (the
only kind of Pad\'e approximant considered in this paper)
are real.

\begin{definition}
The \emph{Pick class} consists of all holomorphic
functions $\UHP \to \overline{\UHP}$.
\end{definition}

Non-constant members of the Pick class map the complex upper half-plane
$\UHP$ into itself, as follows from the open mapping theorem.
The one-to-one correspondence between the Schur and Pick
classes can be used to transform properties of the Schur class,
such as the characterisation of Schur rational maps as finite
Blaschke products multiplied by numbers in $\overline{\D}$,
into statements about the Pick class.
But for our purposes it is simpler to work
directly with the Pick class.

\begin{lemma} \label{lposd}
If a rational map $\mathcal R$ is in the Pick class and
$\mathcal R$ is real-valued (hence finite) at some point
$x \in \R$ then either $\mathcal R$ is a constant or
$D{\mathcal R}(x) > 0$.
\end{lemma}
\begin{proof}
If $D{\mathcal R}(x) < 0$ then all points
$x+\varepsilon i \in \UHP$ with $\varepsilon > 0$
sufficiently small would be mapped into the
lower half-plane.
If $D{\mathcal R}(x) = 0$ and
$\mathcal R$ is not a constant, then
${\mathcal R}(z) =
{\mathcal R}(x) + \alpha (z-x)^k
+ \BigO((z-x)^{k+1})$
with $k > 1$ and $\alpha \neq 0$,
so again some points in $\UHP$
would be mapped into the lower half-plane.
\end{proof}

As the following lemma shows, every real rational Pick map of
positive degree can be generated via the inverse Pick
algorithm from a real rational Pick map of degree one smaller
($\mu$ should be taken positive in the inverse algorithm
in order to generate a map with positive derivative at $x$).

\begin{lemma}[Degree reduction] \label{lrecurse}
Let $\mathcal R$ be a real rational map and $x \in \R$
some point for which ${\mathcal R}(x)$ is finite and
$D{\mathcal R}(x) > 0$.
Then $\mathcal R$ is in the Pick class if and only
if $\Pick{x}{\mathcal R}$ is in the Pick class.
\end{lemma}
\begin{proof}
Let $\mathcal T = \Pick{x}{\mathcal R}$ and
recall the relationship
\begin{equation}
{\mathcal R}(z) = {\mathcal R}(x) + \frac{D{\mathcal R}(x)(z-x)}{1-(z-x){\mathcal T}(z)}. \label{cfrac}
\end{equation}

First suppose that $\mathcal T$ is in the Pick class and
take an arbitrary point $z_0 \in \UHP$.  In order to see that
${\mathcal R}(z_0) \in \UHP$, let $B = {\mathcal T}(z_0) \in \overline{\UHP}$
and consider the M\"obius transformation
\begin{equation}
M_B : z \mapsto {\mathcal R}(x) + \frac{D{\mathcal R}(x)(z-x)}{1 - (z-x) B}.
\end{equation}
Since $z \mapsto -1/z$ maps $\UHP$ to $\UHP$, so does
$z \mapsto B - 1/z$ because $\Im(B) \geq 0$.
Composing with $z \mapsto - D{\mathcal R}(x)/z$
shows that $z \mapsto z\,D{\mathcal R}(x)/(1 - B z)$ also maps
$\UHP$ into itself since $D{\mathcal R}(x) > 0$.  Thus $M_B$ maps $\UHP$ into itself
because $x$ and ${\mathcal R}(x)$ are real.  In particular ${\mathcal R}(z_0) =
M_B(z_0) \in \UHP$.  This shows that ${\mathcal R}$ is in the Pick class.

Now suppose that $\mathcal R$ is in the Pick class.
The Poincar\'e distance $d(z,w)$ between points $z,w$ of $\UHP$ is
given by
\begin{equation}
d(z,w) = \log \frac{|z-\overline{w}| + |z-w|}{|z-\overline{w}| - |z-w|}.
\end{equation}
Since $\mathcal R$ maps $\UHP$ holomorphically into itself, it does not
expand the Poincar\'e distance:
$d({\mathcal R}(z),{\mathcal R}(w)) \leq d(z,w)$, which is equivalent to
\begin{equation} \label{epm}
\frac{|{\mathcal R}(z)-\overline{{\mathcal R}(w)}| + |{\mathcal R}(z)-{\mathcal R}(w)|}{|z-\overline{w}| + |z-w|}
\leq
\frac{|{\mathcal R}(z)-\overline{{\mathcal R}(w)}| - |{\mathcal R}(z)-{\mathcal R}(w)|}{|z-\overline{w}| - |z-w|}.
\end{equation}
Writing $w = x + \varepsilon i$ and passing to the limit $\varepsilon \downarrow 0$
in (\ref{epm}) yields
\begin{equation}
\frac{|{\mathcal R}(z)-{\mathcal R}(x)|}{|z-x|} \leq D{\mathcal R}(x) \frac{\Im({\mathcal R}(z))}{\Im(z)}
\frac{|z-x|}{|{\mathcal R}(z)-{\mathcal R}(x)|}
\end{equation}
(recall that $D{\mathcal R}(x) > 0$).
Using the identity $\Im(\alpha)/|\alpha|^2 = - \Im(1/\alpha)$, this
rearranges to
\begin{equation}
0 \leq \Im\left(\frac{1}{z-x} - \frac{D{\mathcal R}(x)}{{\mathcal R}(z) - {\mathcal R}(x)}\right)
\end{equation}
which is exactly $\Im({\mathcal T}(z)) \geq 0$ since
\begin{equation}
{\mathcal T}(z) = \frac{1 - \frac{D{\mathcal R}(x)(z-x)}
{{\mathcal R}(z) - {\mathcal R}(x)}}{z-x}.
\end{equation}
\end{proof}

It is now easy to understand why a real rational map in the Pick class
has non-negative Schwarzian derivatives of all orders: applying
the Pick algorithm repeatedly gives a sequence of real rational
Pick maps of decreasing degree, finishing with a constant.  Except
for the constant, these all have positive derivative (lemma~\ref{lposd}).
But the higher order Schwarzians of the original map are
just products of these derivatives, up to a positive constant
(corollary~\ref{schsch}).

\begin{lemma}[Characterization] \label{lchar}
Let $\mathcal R$ be a real rational map of degree $d \geq 1$,
and $x \in \R$ some point at which $\mathcal R$ is finite.
If $\mathcal R$ is in the Pick class then $D{\mathcal R}(x) > 0$
and $S_k({\mathcal R})(x) > 0$ for $1 \leq k < d$.
The existence of the Schwarzian derivatives
is part of the conclusion.
Conversely, if $D{\mathcal R}(x) > 0$ and
$S_k({\mathcal R})(x) \geq 0$ for $1 \leq k < d$
then $\mathcal R$ is in the Pick class.
The existence of the Schwarzian derivatives
is part of the hypotheses.
\end{lemma}
\begin{proof}
By induction on the degree.  The case $d=1$ is easily checked,
so take $d > 1$.
If $\mathcal R$ is in the Pick class then
$D{\mathcal R}(x) > 0$ (lemma~\ref{lposd})
and ${\mathcal T} \equiv \Pick{x}{\mathcal R}$
is in the Pick class (lemma~\ref{lrecurse}).
Because $\degree {\mathcal T} = d-1 \geq 1$, it follows
by induction that $D{\mathcal T}(x) > 0$ and
$S_k({\mathcal T})(x) > 0$ for $1 \leq k < d-1$.
Corollary~\ref{schsch}
immediately gives $S_k({\mathcal R})(x) > 0$ for $1 \leq k < d$.

Conversely, if $D{\mathcal R}(x) > 0$ and
$S_k({\mathcal R})(x) \geq 0$ for $1 \leq k < d$,
then, by corollary~\ref{schsch},
$S_k({\mathcal T})(x) \geq 0$ for $1 \leq k < d-1$.  Note that
$S_1({\mathcal R})(x) \neq 0$ (otherwise $\mathcal R$ would not
be normal of order $2$; by hypothesis the Pad\'e
approximants to $\mathcal R$ at $x$ of order $1, \ldots, d$
exist, so corollary~\ref{cnormal} would then imply
that $\mathcal R$ equals $\Pade{\mathcal R}{x}{1}$,
which has degree $1$, a contradiction with $d > 1$),
so $D{\mathcal T}(x) > 0$ by corollary~\ref{schsch}.
Thus by induction $\mathcal T$ is in the Pick class,
and therefore also $\mathcal R$ by lemma~\ref{lrecurse}.
\end{proof}

\begin{corollary} \label{picknorm}
Let $\mathcal R$ be a real rational Pick map of
degree $d$ which is finite at $x$.
Then $\mathcal R$ is normal of orders $0, \ldots, d$ at $x$.
\end{corollary}
\begin{proof}
If $\mathcal R$ is constant then there is nothing to prove.
Otherwise $D{\mathcal R}(x) \neq 0$ by lemma~\ref{lposd},
which means that $\mathcal R$ is normal of order $1$ at $x$.
Since $S_1({\mathcal R})(x) \neq 0$ by lemma~\ref{lchar},
it follows from lemma~\ref{lform} that $\det M_2(x,{\mathcal R}) \neq 0$,
which shows that $\mathcal R$ is normal of order $2$ at $x$
(corollary~\ref{cnorm}).  Repeat for higher orders using
the non-zero Schwarzian derivatives ensured by lemma~\ref{lchar}.
\end{proof}

Thus every real rational Pick map of degree $d$ which is
finite at a point $x \in \R$ can be written in the form
\begin{equation} \label{confrac2}
{\mathcal R}(z) =
A_0 + \cfrac{\mu_0(z-x)}{1 - (z-x) A_1 - \cfrac{\mu_1(z-x)^2}{1 - (z-x) A_2 - \cfrac{\mu_2(z-x)^2}{\quad\cfrac[l]{\ddots}{1 - (z-x)A_d}}}}
\end{equation}
where $A_0, \ldots, A_d \in \R$ and
$\mu_0, \ldots, \mu_{d-1}$ are strictly positive.
Conversely, every function of this form defines a real
rational Pick map of degree $d$ which is finite at $x$.

\begin{corollary} \label{cschpick}
Let $d$ be a positive integer and $x$ some point.
Suppose $f$ has $2d$ derivatives at $x$ and $Df(x) > 0$.
Then $S_1(f)(x) > 0$, $\ldots$, $S_{d-1}(f)(x) > 0$
if and only if the $d$'th Pad\'e approximant to
$f$ at $x$ exists, has degree $d$, and is in the Pick class.
\end{corollary}
\begin{proof}
If $S_1(f)(x) > 0$, $\ldots$, $S_{d-1}(f)(x) > 0$ then,
like in the proof of corollary~\ref{picknorm}, $f$
is normal of orders $0, \ldots, d$ by induction
on the order.  In particular $\Pade{f}{x}{d}$ exists
and has degree $d$.
Now suppose that ${\mathcal R} \equiv \Pade{f}{x}{d}$ exists
and has degree $d$.
Since $f$ and $\mathcal R$ coincide to order $2d$,
it is immediate that $S_k({\mathcal R})(x) = S_k(f)(x)$ for
$1 \leq k < d$.
Thus $S_1(f)(x) > 0, \ldots, S_{d-1}(f)(x) > 0$ if and
only if $S_1({\mathcal R})(x) > 0, \ldots, S_{d-1}({\mathcal R})(x) > 0$,
and this, according to lemma~\ref{lchar}, if and only if ${\mathcal R}$
is in the Pick class.
\end{proof}

\section{Composition formula}
\label{sec:composition-formula}

The composition formula for the classical Schwarzian derivative,
$S_1(g\circ f) = S_1(g)\circ f \, (Df)^2 + S_1(f)$,
implies that the set of maps for which $S_1$ is identically
zero is closed under composition.
This set is precisely the group of M\"obius transformations ---
the composition formula implicitly contains the group
structure of these maps.  On the other hand, the set of
rational maps of degree at most $d$ is not closed
under composition, yet these are the maps for which
$S_d$ is identically zero (lemma~\ref{lzero}).
Inevitably the composition formula for
$S_d$ contains additional terms reflecting the
lack of group structure:

\begin{lemma} \label{lcomp}
Let $d$ be a positive integer, $f$ (resp.\ $g$) a function which is
$2d+1$ times differentiable at $x$ (resp.\ $f(x)$).
Suppose $Df(x) \neq 0$, $Dg(f(x)) \neq 0$ and
$S_d(f)(x)$, $S_d(g)(f(x))$ and $S_d(g \circ f)(x)$
exist.  Then
\begin{equation} \label{ecomp}
S_d(g \circ f)(x) =
S_d(g)(f(x)) (Df(x))^{2d} + S_d(f)(x)
+ S_d(\Pade{g}{f(x)}{d} \circ \Pade{f}{x}{d})(x).
\end{equation}
\end{lemma}

Note that if either $\Pade{g}{f(x)}{d}$ or $\Pade{f}{x}{d}$ is a M{\"o}bius
transformation, then $S_d(\Pade{g}{f(x)}{d} \circ \Pade{f}{x}{d})(x) = 0$
because $\Pade{g}{f(x)}{d} \circ \Pade{f}{x}{d}$ has degree at most $d$.

\begin{proof}
The result is almost immediate from the definitions.
Indeed, set ${\mathcal F} \equiv \Pade{f}{x}{d}$
and ${\mathcal G} \equiv \Pade{g}{f(x)}{d}$;
define $\Delta_f = {\mathcal F}^{-1} \circ f$,
$\Delta_g = {\mathcal G}^{-1} \circ g$
and
$\Delta_\circ = (\Pade{{\mathcal G} \circ {\mathcal F}}{x}{d})^{-1} \circ
{\mathcal G} \circ {\mathcal F}$.
By definition
\begin{eqnarray}
\Delta_f(z) &=& z + \frac{S_d(f)(x)}{(2d+1)!}\,(z-x)^{2d+1} + \littleo((z-x)^{2d+1}), \nonumber \\
\Delta_g(z) &=& z + \frac{S_d(g)(f(x))}{(2d+1)!}\,(z-f(x))^{2d+1} + \littleo((z-f(x))^{2d+1}), \nonumber \\
\Delta_\circ(z) &=& z + \frac{S_d({\mathcal G} \circ {\mathcal F})(x)}{(2d+1)!}\,(z-x)^{2d+1} +
\littleo((z-x)^{2d+1}).
\end{eqnarray}
Let $\widehat{\Delta_g} =
{\mathcal F}^{-1} \circ \Delta_g \circ {\mathcal F}$.
Manipulating the $\Delta_g$ series and using $D{\mathcal F}(x) = Df(x)$,
\begin{equation}
\widehat{\Delta_g}(z) = z + \frac{S_d(g)(f(x))}{(2d+1)!} Df(x)^{2d} \,(z-x)^{2d+1} + \littleo((z-x)^{2d+1}).
\end{equation}
Observe that
$g \circ f = {\mathcal G} \circ \Delta_g \circ {\mathcal F} \circ \Delta_f
= {\mathcal G} \circ {\mathcal F} \circ \widehat{\Delta_g} \circ \Delta_f
= \Pade{{\mathcal G} \circ {\mathcal F}}{x}{d} \circ \Delta_\circ \circ \widehat{\Delta_g} \circ \Delta_f$.
Composing the series for last three terms gives
\begin{equation*}
\Delta_\circ \circ \widehat{\Delta_g} \circ \Delta_f(z) =
\quad\quad\quad\quad\quad\quad\quad\quad\quad\quad
\quad\quad\quad\quad\quad\quad\quad\quad\quad\quad
\quad\quad\quad\quad\quad
\end{equation*}
\begin{equation*}
z + \left(S_d({\mathcal G} \circ {\mathcal F})(x)+S_d(g)(f(x))Df(x)^{2d} + S_d(f)(x)\right)
\,\frac{(z-x)^{2d+1}}{(2d+1)!}
\end{equation*}
\begin{equation*}
\quad\quad\quad\quad\quad\quad\quad\quad\quad\quad
\quad\quad\quad\quad\quad\quad\quad\quad\quad\quad
\quad\quad\quad\quad\quad
+ \littleo((z-x)^{2d+1}).
\end{equation*}
Thus
$\Pade{g\circ f}{x}{d} = \Pade{{\mathcal G} \circ {\mathcal F}}{x}{d}$
and equation~\ref{ecomp} follows.
\end{proof}

Nonetheless, for any $d$ there is a general composition inequality for
maps with non-negative Schwarzian derivatives of lower order.  The reason
for this is that Pad\'e approximants to such maps correspond to (non-constant)
members of the Pick class, and the set of such Pick maps is closed under
composition.

\begin{proposition} \label{pcomp}
Let $d$ be a positive integer and $f$ (resp.\ $g$) a function which is
$2d+1$ times differentiable at $x$ (resp.\ $f(x)$).
Suppose $Df(x) \neq 0$, $Dg(f(x)) \neq 0$,
$S_d(f)(x)$ and $S_d(g)(f(x))$ exist, and
$S_k(f)(x) \geq 0$, $S_k(g)(f(x)) \geq 0$ for
every $1 \leq k < d$.  Then $S_d(g \circ f)(x)$
exists and
\begin{equation}
S_d(g \circ f)(x) \geq S_d(g)(f(x)) (Df(x))^{2d} + S_d(f)(x).
\end{equation}
\end{proposition}
\begin{proof}
Without loss of generality $Df(x) > 0$ and $Dg(f(x)) > 0$
(otherwise pre- and/or post-compose with $z \mapsto -z$
to arrange this --- Schwarzian derivatives do not change).
The $d$'th Pad\'e approximant ${\mathcal F} \equiv
\Pade{f}{x}{d}$ to $f$ at $x$ is in the Pick class
by lemma~\ref{lchar} because $S_k({\mathcal F})(x) =
S_k(f)(x) \geq 0$ for $1 \leq k < d$.  Likewise
${\mathcal G} \equiv \Pade{g}{f(x)}{d}$ is in the Pick
class.
In order words, $\mathcal F$ and $\mathcal G$
map the complex upper half-plane $\UHP$ into itself.
Obviously their composition does too, i.e.\ ${\mathcal G}
\circ {\mathcal F}$ is in the Pick class.
Thus $\Pade{{\mathcal G} \circ {\mathcal F}}{x}{d}$
exists by lemma~\ref{picknorm} --- this is equivalent
to the existence of $S_d(g \circ f)$.  Finally,
$S_d({\mathcal G} \circ {\mathcal F})(x) \geq 0$ by
lemma~\ref{lchar}.
\end{proof}

\section{Monotone matrix functions} \label{smonmat}

In this section we introduce the class of monotone matrix functions
and relate them to maps with non-negative higher order Schwarzian
derivatives.  A complete description of this class can be
found in~\cite{DO}.

Recall how to take the image of a real symmetric
matrix $A$ by a function $f$: if $A$ is diagonal,
$A = \mathrm{diag}(\lambda_1, \ldots, \lambda_n)$,
then $f(A) = \mathrm{diag}(f(\lambda_1), \ldots, f(\lambda_n))$;
otherwise diagonalise $A$ via some linear coordinate change,
take the image of the diagonal matrix, and undiagonalise
by applying the inverse coordinate change.  This is
well-defined if the spectrum of $A$ belongs to the
domain of $f$.

Recall the ordering on the real symmetric $n$-by-$n$ matrices:
$A \leq B$ if and only if $B - A$ is a positive matrix, meaning
$v^T (B-A) v \geq 0$ for every $n$-by-$1$ vector $v$.

\begin{definition}
Let $n$ be a positive integer, $U$ an open interval of the
real line, and $f : U \to \R$ a function.
Call $f$ \emph{matrix monotone of order $n$}
if, for any real symmetric $n$-by-$n$ matrices $A$ and $B$ with
spectrum contained in $U$,
$A \leq B$ implies $f(A) \leq f(B)$.
\end{definition}

\begin{lemma} \label{lmm}
Let $d$ be a positive integer, $U$ an open interval
and $f : U \to \R$ a function with $2d+1$ derivatives.
Suppose $Df(x) > 0$ for all $x \in U$.
Then $S_k(f) \geq 0$ on $U$ for all $1 \leq k \leq d$
if and only if $f$ is matrix monotone of order $d+1$.
\end{lemma}
\begin{proof}
Suppose first that $S_k(f) \geq 0$ for all $1 \leq k \leq d$.
Take some $x \in U$ and let $\mathcal R$ be the $d$'th Pad\'e
approximant to $f$ at $x$.  This rational map is in the Pick
class because $S_k({\mathcal R})(x) = S_k(f)(x) \geq 0$
for $1 \leq k < d$ (lemma~\ref{lchar}),
so $M_{d+1}(x,{\mathcal R})$ is a positive matrix
by~\cite[theorem III.IV]{DO}.  Since $f$ and $\mathcal R$
coincide to order $2d$ at $x$ and
$D^{2d+1}(f)(x) = D^{2d+1}({\mathcal R})(x) +
Df(x) S_d(f)(x)$ (this is the alternative definition
of the Schwarzian derivative from section~\ref{ratapprox}),
\begin{equation}
M_{d+1}(x,f) = M_{d+1}(x,{\mathcal R}) +
Df(x) S_d(f)(x)
\left[
\begin{array}{ccc}
0 & \cdots & 0 \\
\vdots && \vdots \\
0 & \cdots & 1
\end{array}
\right]
\geq M_{d+1}(x,{\mathcal R}).
\end{equation}
Thus $M_{d+1}(x,f)$ is also a positive matrix.
Since $M_{d+1}(x,f)$ is positive for every $x \in U$,
theorem~VIII.V of~\cite{DO} shows%
\footnote{The convexity hypothesis in the theorem is used to
get the existence of sufficiently many derivatives almost
everywhere; if the function is assumed sufficiently differentiable,
as here, then this hypothesis is automatically satisfied --- it
follows from the positivity of the matrix.}
that $f$ is matrix monotone of order $d+1$.

Now suppose that $f$ is matrix monotone of order $d+1$ on $V$.
Choose some point $x \in U$.  Applying~\cite[theorem XIV.I]{DO}
with $n=d+1$ and $S$ consisting of $2d+1$ copies of $x$ gives a
Pick function $\phi$ on $U$ which coincides with $f$ to order $2d$ at
$x$.  According to theorem~III.IV of~\cite{DO}, either
$M_{d}(x,\phi)$ is strictly positive, or $\phi$ is a
rational map of degree at most $d-1$.
In this last case, \cite[theorem XIV.II]{DO} states that
$f$ and the rational Pick map $\phi$ coincide on $U$;
the result is then immediate from lemma~\ref{lchar}.
So suppose $M_{d}(x,\phi)=M_d(x,f)$ is strictly positive.
Then the principal minors of $M_d(x,f)$ are strictly positive:
$\det M_j(x,f) > 0$ for $1 \leq j \leq d$.
Thus $f$ is normal of orders $1, \ldots, d$
at $x$ and formula~\ref{esderiv} can be freely applied.
This gives $S_k(f)(x) > 0$ for $1 \leq k < d$.
Furthermore, $M_{d+1}(x,f)$ is a positive matrix
by~\cite[theorem VII.VI]{DO}, so $S_d(f)(x) \geq 0$.
\end{proof}

It is a remarkable fact (Loewner's theorem~\cite{DO}) that a
function is matrix monotone of all orders if and only if it
extends holomorphically to the complex upper half-plane $\UHP$
and maps $\UHP$ into itself.
Thus:
\begin{proof}[Proof of proposition~\ref{puni}]
Combine lemma~\ref{lmm} and Loewner's theorem.
\end{proof}

\section{Proof of the generalised Koebe lemma} \label{skoebe}
In this section we prove theorem~\ref{koebe}.  Without loss of
generality $Df > 0$ on $U$, so $f$ is matrix monotone of order
$d+1$ (lemma~\ref{lmm}).
Applying~\cite[theorem XIV.I]{DO} with $n=d+1$ and $S$ consisting
of $2d+1$ copies of $x$ gives a Pick function $\phi$ on $U$
which coincides with $f$ to order $2d$ at $x$.  It clearly suffices
to prove the result for $\phi$.  Now $\phi$, being in the Pick class,
has the integral representation
\begin{equation}
\phi(x)=\alpha x + \beta + \int \left[ \frac 1{\xi-x} - \frac
  \xi{\xi^2+1}\right]\,d\mu(\xi),
\end{equation}
where $\alpha \ge 0$, $\beta$ is real and $\mu$ is a
positive Borel measure, supported in $\R \setminus U$,
for which $\int (\xi^2+1)^{-1}\,d\mu(\xi)$ is finite.
See~\cite[theorem II.I, lemma II.2]{DO}.
The derivatives of $\phi$ are given by the following formulae:
\begin{eqnarray}
D\phi(x)   &=& \alpha + \int \frac 1{(\xi-x)^2}\,d\mu(\xi) \nonumber \\
D^m\phi(x) &=& m! \int \frac 1{(\xi-x)^{m+1}}\,d\mu(\xi) \,\,\,\mbox{
  if $m>1$.}
\end{eqnarray}
We can now estimate $|D^m \phi(x)|$ easily:
\begin{eqnarray}
  |D^m \phi(x)| &=& m! \left|\int \frac 1{(\xi-x)^{m+1}}\,d\mu(\xi)\right| \nonumber\\
  &\le& m! \int \frac 1{|\xi-x|^{n+1}} \frac 1{|\xi -x|^{m-n}}\,d\mu(\xi) \nonumber\\
  &\le& \frac {m!}{n!} \frac 1{\dist(x, \partial U)^{m-n}} |D^n \phi(x)|.
\end{eqnarray}
The last assertion holds because $n$ is odd, so $|\xi-x|^{n+1}=(\xi
-x)^{n+1}$, and because $\mu$ puts no mass on the interval $U$. Also
notice that in the case $n=1$ we have used the positivity of $\alpha$.

\section{Proof of the main theorem} \label{smain}

In this section we prove theorem~\ref{tmain}.  Let us recall a few
definitions. An interval is called \emph{nice} if the iterates of its
boundary points never return inside the interval.  A sequence of
intervals $\{W_j\}_{j=0}^s$ is called a {\it chain} if $W_{j}$ is a
connected component of $f^{-1}(W_{j+1})$.
The {\it order} of a chain is the number of intervals in the chain
containing a critical point.
A $\delta$-scaled neighbourhood of an interval $J$ is any $V$
containing the set $\{x : \exists y \in J, |x-y| < \delta |J|\}$.
In this case we also say that $J$ is $\delta$-well-inside $V$.
In what follows we will assume that $I$ is the interval $[0,1]$.

We need real bounds, and will use theorem~D$'$
of \cite[section~8]{SV}.  The authors formulate this result slightly
differently, however their proof gives precisely this:

\begin{fact} \label{f1} Let $f:I\to I$ be a $C^3$ map with non flat
  critical points. Then there exists $\tau>0$ and arbitrarily small
  neighbourhoods $\W$ and $\V \supset \W$ of the set of those critical
  points which are not in the basin of any periodic attractor such
  that
  \begin{itemize}
  \item all connected components of $\W$ and $\V$ are nice intervals;
  \item if $W\subset V$ are two connected components of $\W$ and $\V$,
    then $V$ is a $\tau$-scaled neighbourhood of $W$;
  \item if $x \in I$ is a point and $s \geq 0$ is minimal such that
    $f^s(x)\in \W$, then the chain obtained by pulling back the connected
    component of $\V$ containing $f^s(x)$ along the orbit
    $x,f(x),\ldots, f^s(x)$ has order bounded by the number of critical
    points of $f$.
  \end{itemize}
\end{fact}

We also make use of theorem~2 from \cite{ST}, which we
state as:

\begin{fact} \label{sht}
  Let $f:I\to I$ be a $C^n$ map with non flat critical points, $n\ge 2$.
  Let $T$ be an interval
  such that $f^s : T \to f^s(T)$ is a diffeomorphism. For each
  $S,\tau,\epsilon > 0$  there exists $\delta =
  \delta(S, \tau, \epsilon, f) > 0$ 
  satisfying the following. If $\sum_{j=0}^{s-1} |f^j(T)| \le S$ and
  $J$ is a subinterval of $T$ such that
  \begin{itemize}
  \item
    $f^s(T)$ is a $\tau$-scaled neighbourhood of $f^s(J)$;
  \item ~$|f^j(J)| < \delta$  for $0\le j < s$,
  \end{itemize}
  then, letting $\phi_0 : J \to I$ and $\phi_s : f^s(J) \to I$ be affine
  diffeomorphisms, there exists a real-analytic diffeomorphism $G : I \to I$
  such that $\|\phi_s f^s \phi^{-1}_0 - G\|_{C^n} <\epsilon$,
  and $G^{-1}$ belongs to the $\PoSch_\infty((-\tau/2,1+\tau/2))$ class.
\end{fact}

In addition, we will use the following lemmas:

\begin{lemma} \label{icmp}
  Let $F,G : I\to I$ be two $C^n$ diffeomorphisms, $n \geq 1$,
  with $\|F - G\|_{C^n} <\epsilon$.
  Take some $K, \epsilon > 0$ and suppose
  $|DG(x)|> K^{-1}$ and $|D^kG(x)| < K$
  for all $x \in I$ and $k=1,\ldots,n$.
  Then there exists $\delta=\delta(n, K, \epsilon)$ such that
  $\|F\circ G^{-1} - Id\|_{C^n} < \delta$.
  Moreover, $\lim_{\epsilon \to 0} \delta (n, K, \epsilon)=0$.
\end{lemma}

\begin{proof}
Denote $H=F\circ G^{-1}$.
Then $|H(x) - x| = |F - G|(G^{-1}(x)) < \epsilon$.
The first derivative of $H$ is
$DH(x)= (DF/DG)(G^{-1}(x))$
so obviously $|DH(x) -1|< \epsilon K$.
This proves the lemma for $n=1$ with $\delta(1, K, \epsilon) = \epsilon K$.
For $n \geq 2$ we reason inductively.  Observe that
$$
(D^nH) \circ G = \quad\quad\quad\quad\quad\quad\quad\quad\quad\quad\quad\quad\quad\quad\quad\quad\quad\quad\quad\quad\quad\quad\quad\quad\quad\quad\quad
$$
$$
\quad\quad\frac{D^nF - D^nG + D^nG \left(1 - H'\circ G\right) -
\sum_{i=2}^{n-1} (D^iH)\circ G Q_{n,i}(DG, \ldots, D^{n+1-i}G)}
{(DG)^n}
$$
where the $Q_{n,i}$ are polynomials (this is easily checked by induction).

Since $\|DG\|_{C^{n-1}} < K$ by hypothesis, there is some
$\Delta = \Delta(n, K)$ for which
$|Q_{n,i}(DG, \ldots, D^{n+1-i}G)| \leq \Delta$
for $i = 2, \ldots, n-1$.
Thus
\begin{eqnarray}
|D^nH| & \leq &
\frac{\|F-G\|_{C^n} + \|DG\|_{C^{n-1}}\, \epsilon K + \sum_{i=2}^{n-1} \|H - Id\|_{C^i}\Delta(n, K)}{K^{-n}}
\nonumber \\
& \leq & K^n\left(\epsilon + \epsilon K^2 + n\,\delta(n-1, K, \epsilon)\,\Delta(n, K) \right)
\equiv \delta'(n, K, \epsilon).  \label{CnBound}
\end{eqnarray}

Clearly $\delta'(n, K, \epsilon) \to 0$ as $\epsilon \to 0$.
The result follows.
\end{proof}

\begin{lemma} \label{pert}
For every $A \geq 0$, $\alpha > 1$ and $d \geq 1$ there is some
$\epsilon = \epsilon(A, \alpha, d) > 0$ with the following property.
If $\phi, \psi : I \to I$ are $C^{2d+1}(I)$
diffeomorphisms,
$\|\phi - Id\|_{C^{2d+1}} <\epsilon$,
$\|\psi - Id\|_{C^{2d+1}} <\epsilon$
and $0 \leq a \leq A$,
then $(\psi \circ q_{\alpha,a} \circ \phi)^{-1} \in \PoSch_d(\interior I)$ where
\begin{equation}
q_{\alpha,a} : x \mapsto
\frac{(x+a)^\alpha - a^\alpha}{(1+a)^\alpha - a^\alpha}.
\end{equation}
\end{lemma}
Recall that we have assumed $I=[0,1]$.

\begin{proof}
The only difficulty here is that $q_{\alpha,a}$ has a singularity
at $-a$ which may be arbitrarily close to $I$.  This means that
derivatives of $q_{\alpha,a}^{-1}$ are not uniformly bounded as $a \to 0$,
complicating continuity arguments.  In what follows we only mark dependence
on $a$ explicitly: the other parameters $\alpha$,
$A$ and $n$ should be considered as fixed, with all quantities
potentially depending on them.  For example, we will write
$q_a$ for $q_{\alpha,a}$.

The singularity can be side-stepped by decomposing $q_a$
as $t \circ s_a \circ r$ where $s_a : I \to I$ is a
real-analytic homeomorphism, diffeomorphic on $\interior I$,
with $s_a^{-1}$ in the Pick class on $\interior I$, like $q_a$.
The maps $r, t : I \to I$ should be real-analytic diffeomorphisms,
with $r^{-1}$ and $t^{-1}$ Pick functions on $J$,
an interval strictly bigger than $I$.  Finally, $r^{-1}$ and
$t^{-1}$ should not be rational functions.  Then the
the matrices $M_k(y,r^{-1})$ and $M_k(y,s^{-1})$ will be uniformly
positive for $y \in I$ and $k = 1, \ldots, n$
thanks to~\cite[theorem~III.IV]{DO}.
Note that $r$, $t$ and $J$ do not depend on $a$.
Although it is not hard to give an explicit such decomposition,
we will not do so here since the formulae are ugly and
uninformative.

With such a decomposition in hand, \cite[theorem~VII.V]{DO}
and an easy continuity argument imply that
$(\psi \circ t)^{-1} \in \PoSch_d(I)$ and $(r \circ \phi)^{-1} \in \PoSch_d(I)$
if $\|\phi - Id\|_{C^{2d+1}}$ and $\|\psi - Id\|_{C^{2d+1}}$ are
sufficiently small, which is the desired result.
\end{proof}

Now we can finish the proof of the theorem.
The notation $g \in \PoSch_d^{-1}(x)$ means
$S_1(g^{-1})(g(x)) \geq 0, \ldots, S_d(g^{-1})(g(x)) \geq 0$
(local inverse at $x$).
Assume that a critical point $c$ is not contained in the basin of a periodic
attractor. The neighbourhood $X$ of $c$ will be a connected component of $\W$ given
by fact~\ref{f1} for $\W$ sufficiently small.

Take $\W$ and $\V$ as in fact~\ref{f1} and let $f^s(x) \in \W$.  We
may suppose that $s \geq 0$ is minimal with $f^s(x) \in \W$, since the
general case can be deduced from this by decomposition (recall that
the class $\PoSch_d$ is closed under composition).

Let $V$ and $W$ be connected components of $\V$ and $\W$ containing $c$.
Let $V'$ be a $\tau/3$--scaled neighbourhood of $W$ so that
$|V'|=(1+2/3\tau)|W|$.  The interval $V$ is also a $\tau/3$--scaled
neighbourhood of $V'$.
Let $\{V_j\}_{j=0}^s$ be the corresponding chain of pullbacks of $V$
along the orbit of $x$, i.e. $f^j(x)\in V_j$ and $V_s=V$, and let
$\{V'_j\}_{j=0}^s$ and $\{W_j\}_{j=0}^s$ be corresponding chains for $V'$
and $W$.
The map $f^s : W_0 \to W$ is a diffeomorphism because of the
minimality of $s$.  This is true even if $\W$ does not contain every
critical point, as long as $\W$ is sufficiently small (the necessary
smallness does not depend on $x$ or $s$).

Due to \cite[theorem C]{SV} there exist $\tau'>0$ and $C>0$ such that
$V'_j$ is a $\tau'$--scaled neighbourhood of $W_j$, $V_j$ is a
$\tau'$--scaled neighbourhood of $V'_j$ and $|V'_j| < C |W_j|$ for all
$j=0,\ldots,s$. Minimality of $s$ implies that all intervals $W_j$ are
disjoint, so
$$ \sum_{j=0}^s |V'_j| < C$$
and we can use fact~\ref{sht} for $f^k: V'_j\to V'_{j+k}$ if this map
is a diffeomorphism.

Let $0 \leq s_0 < s_1 < \cdots < s_k = s$ be the moments $j$ when
$V'_j$ contains a critical point, and put $s_{-1} = -1$ for
convenience.  Then $f^{s_i - s_{i-1} - 1} : V'_{s_{i-1}+1} \to V'_{s_i}$
is a diffeomorphism for $i = 0, \ldots k$.
It is enough to show that each $f^{s_i - s_{i-1}}$ belongs to
$\PoSch_d^{-1}(f^{s_{i-1}+1}(x))$ for $i = 0, \ldots k$,
since these compose to give $f^{s+1}$.

Take some $0 \leq i < k$ (see below for the case $i=k$)
and apply fact~\ref{sht} to
$f^{s_i - s_{i-1} - 1} : V'_{s_{i-1}+1} \to V'_{s_i}$
with $T = V'_{s_{i-1}+1}$ and $J = W_{s_{i-1}+1}$.
Let $F_i : I \to I$ be $f^{s_i - s_{i-1} - 1}$
pre- and post-composed by affine maps taking $I$ to
$W_{s_{i-1}+1}$ and $W_{s_i}$ to $I$ respectively.
We obtain a diffeomorphism $G_i : I \to I$ with $G_i^{-1}$ in
the $\PoSch_\infty ((-\tau/6, 1+\tau/6))$  class such that
$\|F_i - G_i\|_{C^{2d+1}} < \epsilon$.  The complex
Koebe lemma gives the bounds on $G_i$ needed to apply
lemma~\ref{icmp}, yielding
$\|F_i \circ {G_i}^{-1} - Id\|_{C^{2d+1}} < \delta$.
Note that by shrinking the neighbourhood $\W$ we can
make $\epsilon$ and $\delta$ as small as we like.

Abusing the notation, let $c$ be a critical point contained in $V'_{s_i}$.
If $\W$ is small enough, it will be the only critical point.
Note that $c$ is not contained in $W_{s_i}$ because
we are considering the case $s_i \neq s$.
Let $F : I \to I$ be $f : W_{s_i} \to f(W_{s_i})$
pre- and post-composed by affine maps taking
$I$ to $W_{s_i}$ and $f(W_{s_i})$ to $I$ respectively,
as in the previous paragraph.
Then $F$ can be written in the form
$\phi \circ q_{\alpha, a} \circ \psi$,
where $q_{\alpha, a}(x)=((x+a)^\alpha - a^\alpha)/((1+a)^\alpha -
a^\alpha)$,
and $\phi$, $\psi$ are diffeomorphisms of $I$ which are
close to the identity map in the $C^{2d+1}$ topology if
$V'_{s_i}$ is small.
These assertions on $F$ follow from the definition
of a critical point being non-flat.

As noted above, the intervals $V'_{s_i}$ and $W_{s_i}$ are comparable,
so $c$ is not far away from the interval $W_{s_i}$ compared to its size.
Expressed in terms of the rescaled map $F$, this means
that there exists a uniform constant $A>0$ such that the parameter $a$ is
always in $[0,A]$.
Applying lemma~\ref{pert}, we see that the inverse of
$\phi \circ q_{\alpha,a} \circ (\psi \circ F_i \circ G_i^{-1})$
is in $\PoSch_d(I)$, at least if $\W$ is small enough.
Composing with $G_i$, we see that the inverse of
$\phi \circ q_{\alpha,a} \circ \psi \circ F_i$
is also in $\PoSch_d(\interior I)$.  Since this composition is precisely
$f^{s_i-s_{i-1}} : W_{s_{i-1}+1} \to W_{s_i+1}$ rescaled affinely,
this shows that $f^{s_i-s_{i-1}}$ belongs to $\PoSch_d^{-1}(f^{s_{i-1}+1}(x))$
as claimed.

We now consider the case $i=k$, when $s_i=s$.
The difference here is that the critical point
belongs to $W_s = W$, which actually makes the argument slightly simpler.
The critical point cuts the interval $W$ in half,
let $W'$ denote the half containing $f^s(x)$.
Then we repeat the above argument,
but instead of rescaling $f: W \to f(W)$,
we instead rescale $f: W' \to f(W')$;
the parameter $a$ is then always zero.
The argument is otherwise essentially the same.

\end{document}